\documentclass[12pt]{article}
\usepackage{amsmath}
\usepackage{graphicx,psfrag,epsf}
\usepackage{xcolor}
\usepackage{enumerate}
\usepackage{natbib}
\usepackage{url} 

\newcommand{\blind}{0}

\newcommand{\reali}{{\rm I}\negthinspace {\rm R}}

\begin{document}

\bibliographystyle{natbib}

\def\spacingset#1{\renewcommand{\baselinestretch}%
{#1}\small\normalsize} \spacingset{1}


\if0\blind
{
  \title{\bf Likelihood, Replicability and Robbins' Confidence Sequences}
  \author{Luigi Pace\thanks{
    This research was supported by a grant from a Project of National Interest, Italian Ministry of Education, PRIN 2015 -- 2015EASZFS 003, and by the University of Udine grant PRID2017\_DIES004.}\hspace{.2cm}\\
    Department of Economics and Statistics, University of Udine, Italy\\
    and \\
    Alessandra Salvan\thanks{
    This research was supported by a grant from a Project of National Interest, Italian Ministry of Education, PRIN 2015 -- 2015EASZFS 003, and by the University of Padova grant BIRD185955.}\\
    Department of Statistical Sciences, University of Padova, Italy}
 \date{} \maketitle
} \fi

\if1\blind
{
  \bigskip
  \bigskip
  \bigskip
  \begin{center}
    {\LARGE\bf Likelihood,  Replicability and Robbins' Confidence Sequences}
\end{center}
  \medskip
} \fi

\bigskip
\begin{abstract}
The widely {claimed} replicability crisis in science may lead to revised standards of significance. The customary  frequentist confidence intervals,    calibrated through hypothetical repetitions of the experiment that is supposed to have produced the data at hand, rely on a feeble  concept of replicability.   In particular, contradictory  conclusions  may be  reached  when a substantial enlargement of the study is undertaken.   
To redefine statistical confidence in such a way that inferential conclusions are non-contradictory, with large enough probability, under enlargements of the sample,   
 we give a new reading of  a proposal dating back to the 60's, namely Robbins' confidence sequences.  Directly bounding the probability of reaching, in the future, conclusions that contradict the current ones, Robbins' confidence sequences ensure a clear-cut form of replicability when inference is performed on accumulating data. Their main frequentist property is easy  to understand and to prove.  We show that Robbins' confidence sequences may be justified under various views of inference: they are likelihood-based,  can incorporate prior information,  and obey the strong likelihood principle. They are easy to compute, even when inference is on a  parameter of interest, especially using a closed-form approximation from normal asymptotic theory.  
\end{abstract}

\noindent%
{\it Keywords:}  Confidence region; Laplace expansion; Profile likelihood; Revision of standards; Statistical evidence. 
\vfill

\newpage
\spacingset{1.45} 
\section{Introduction}
\label{sec:intro}
Announcing a result is a hazard when the supporting evidence is statistical in nature. In the long run, scientific credibility is undermined if discoveries are claimed  (or understood)  to be more firmly established than they will eventually prove to be. 
The issue {appears to be pressing, especially in the context of the replicability crisis in science claimed by Ioannidis (2005) and many others on its wake. This led to  the  ASA statement,  Wasserstein and Lazar (2016), and to the subsequent 2019 {\it The American Statistician}'s special issue {\it Statistical Inference in the 21st Century: A World Beyond p $<$ 0.05}}.  {General warnings against  misuse and misinterpretation of $p$-values are given also in McShane and Gal (2017) and Kuffner and Walker (2019). Difficulties of objective Bayesian inference in attaining replicability are discussed in Fraser et al.\ (2016).}
 To reduce failure to replicate,   {one solution suggested} in the literature is the use of stricter evidential thresholds, possibly variable by discipline (Johnson, 2013; Goodman, 2016).  Benjamin et al.\ (2017) {and Bickel (2019)} advocate changing the standard threshold for significance from 0.05 to 0.005{, or even 0.001}, while
Lakens et al.\ (2017) recommend a case by case transparently-justified choice, better if   pre-registered.  

 When interest lies in reporting  effect sizes and  related confidence intervals (see e.g.\ Nakagawa and Cuthill, 2007), a revision of standards for statistical significance would entail a parallel revision of standards for confidence levels, say from  0.95 to 0.995. However,  these higher levels are not directly linked to some explicit replicability-related requirement.  {A widely agreed definition of replicability is ``the ability of a researcher to duplicate the results of a prior study if the same procedures are followed but new data are collected'' (Schwalbe, 2016, p.\ 4).  For inference based on confidence regions we introduce a connected, but apparently new,} concept of replicability, 
and explore its relations with a proposal dating back to the 60Õs, namely Robbins' {\it confidence sequences} (Robbins, 1970; see also Darling and Robbins, 1967a,b). 
{To be specific, we require that non-contradictory conclusions are reached when the sample is enlarged,   i.e.\ when information increases. Inferential conclusions from  confidence regions for the same parameter are non-contradictory if these regions overlap; they are instead  contradictory if their intersection is empty. }

Fixed level confidence regions, even with  revised higher levels, fail to fulfil the {non-contradiction}  requirement. 
As a simple example, consider i.i.d.\ sampling from a normal distribution with known variance $\sigma^2_0$.
Let $\bar{Y}_n= \sum_{i=1}^n Y_i/n$ be the sample mean.
Then
$$
\bar{Y}_{n+m} - \bar{Y}_n \sim N \left( 0, \sigma^2_0 \left( \dfrac{1}{n} - \dfrac{1}{n+m} \right) \right)
$$
and the probability that $(1-\alpha)$-level confidence intervals for the mean at sample sizes $n$ and  $n+m$ do not overlap  is
\begin{eqnarray*}
&& 2P_\mu \left( \bar{Y}_{n} +\dfrac{\sigma_0}{\sqrt{n}} z_{1-\alpha/2} < \bar{Y}_{n+m} -\dfrac{\sigma_0}{\sqrt{n+m}} z_{1-\alpha/2} \right)\\
&&= 2 \Phi \left( -z_{1-\alpha/2}\left( \sqrt{1+ \dfrac{n}{m}} +\sqrt{\dfrac{n}{m}} \right) \right) >0\,.
\end{eqnarray*}
Therefore the probability is 1 that we can find  a pair of disjoint intervals,
and consequently it is almost sure that we observe a sequence of samples giving rise to contradictory $(1-\alpha)$-level confidence intervals.
When the realistically attainable sample sizes are very large but finite,
though the usual confidence intervals shrink towards the true value of the parameter as the sample size increases, conflicting conclusions may be reported at various stages of the data acquisition process, with a probability that may be close to 1.

Non-contradiction  is especially compelling in experimental sciences when inference is performed on accumulating data.  Early conclusions are susceptible to be falsified within the matter of years or months, and sometimes even earlier.   When the true state of nature, or a much more reliable representation of it, becomes eventually available,  reputational penalty ensuing from hasty announcement of wrong conclusions could be large. This risk is not present in hard sciences alone.  Think for instance of estimating  the result of an election from early reporting counting areas, where the estimate is made only hours before a winner is declared.    Other contexts where coherence under sample enlargement seems to be cogent are  long-term epidemiological studies and drugs surveillance. {Also on-line randomized experiments (A/B tests) represent a relevant setting}.   

{Even in fixed sample size inference, where sample enlargement is merely hypothetical, non-contradiction may represent a sensible requirement for  replicability-related frequentist evaluation of confidence regions. The requirement is that re-evaluation of a confidence region using an enlarged sample should lead to non-contradiction with a controlled probability.}

In this paper, we show that the use of Robbins' confidence sequences produces non-contradictory confidence regions with probability greater than a fixed lower bound, at least in the idealized situation of i.i.d.\ sampling from a correctly specified parametric model. Robbins' papers are highly technical, and research on confidence sequences seems to have been neglected after the equally technical contributions
Lai (1976) and Csenki (1979).  We try to give an  accessible account  and to highlight the potential impact of Robbins' confidence sequences on principles of statistical inference. We think that their impact should be much larger. Indeed, they can be justified  under various views of inference. These sequences are likelihood-based,  can incorporate prior information, have frequentist properties, have Bayesian properties under a proper prior, and obey the strong likelihood principle. 
Moreover, Robbins' confidence sequences have great pedagogical benefits. They need {no sampling distribution calculations and may require a fairly limited amount of numerical evaluations of the likelihood function}.

 

The outline of the paper is as follows. A new reading of Robbins' confidence sequences is given in Section 2, together with  a closed-form approximation from normal asymptotic theory for a scalar parameter of interest. Inferential  properties of Robbins' confidence sequences are summarized in Section 3, with technical details provided in the Appendix.   Section 4 presents examples  dealing with the normal mean and binomial probabilities and illustrates,  through simulation, the  properties of Robbins' confidence sequences  and  the proposed closed-form approximation. Some conclusions are given in Section 5.  

\section{Non-contradiction and Robbins' confidence sequences}

Let us consider the idealized situation of a statistician who is potentially able to obtain any number, $n$, of observations $y^{(n)}=(y_1,\ldots,y_n)$, realization of the random vector  $Y^{(n)}=(Y_1,\ldots,Y_n)$, not necessarily with i.i.d.\ components. Let $P_\theta$ denote the joint probability distribution of the sequence $Y^{(\infty)}=(Y_1,Y_2,\ldots)$. We suppose that $P_\theta$ belongs to a statistical model with parameter space $\Theta\subseteq \reali^p$.  Let 
$p_n(y^{(n)};\theta )$ denote the density (or probability mass function) of  $Y^{(n)}$ under $P_\theta$. 
Assume that, for every given $n$, all these densities are strictly positive  {on a  support that does not depend on $\theta$}.  

A confidence region, based on $y^{(n)}$ and constructed according to a certain rule,  is a subset of $\Theta$ denoted by $\hat\Theta_n=\hat{\Theta}( y^{(n)})$. A confidence sequence is a sequence of confidence regions. To avoid triviality, we consider only confidence sequences that are consistent, i.e.\ such that $\lim_{n \to \infty} P_\theta ( \theta' \in \hat\Theta_n) =0$ for every $\theta' \neq \theta$, where $\theta, \theta' \in \Theta$.  Consistency implies that, for $\theta'\neq \theta$, 
$$
P_\theta\left( \theta' \in \cap_{n\geq 1} \hat{\Theta}_n  \right) \leq \lim_{n \to \infty} P_\theta ( \theta' \in \hat\Theta_n)= 0\,,
$$
i.e., that a false parameter value cannot belong to all confidence regions  of a consistent confidence sequence. Only the true $\theta$ may belong to $\cap_{n \geq 1}\hat\Theta_n$, {provided  that $\cap_{n \geq 1}\hat\Theta_n$ is non-empty}.


We will say that  a confidence sequence is  non-contradictory, or persistent,  if no $\hat\Theta_n$ is contradicted by a $\hat\Theta_{n+m}$, for some $m>0$. Contradiction happens when, for  an $m>0$,  $\hat\Theta_n \cap \hat\Theta_{n+m}= \emptyset$. When a confidence sequence is non-contradictory, there are conclusions that are common to all confidence statements,  i.e., $\cap_{n \geq 1}\hat\Theta_n \neq \emptyset$.  Consistency ensures that non-contradictory sequences shrink towards the true parameter value.  

Since $\cap_{n \geq 1}\hat\Theta_n = \emptyset$ implies $\theta\notin \cap_{n \geq 1}\hat\Theta_n$, we have
$$
P_\theta\left(  \cap_{n \geq 1}\hat\Theta_n = \emptyset \right) \leq
P_\theta \left(  \theta\notin \cap_{n \geq 1}\hat\Theta_n \right) = 
1- P_\theta \left( \theta \in \hat\Theta_n\;\;{\rm for\;\; every}\;\; n \geq 1  \right)\,.
$$
It follows that, if, for $0 < \varepsilon <1$, 
\begin{equation}\label{persistent}
P_\theta\left( \theta \in \hat{\Theta}_n  \;\;{\rm for\;\; every}\;\; n \geq 1  \right) \geq 1 -\varepsilon  \,,
\end{equation}
then 
$$
P_\theta\left(  \cap_{n \geq 1}\hat\Theta_n = \emptyset \right) \leq \varepsilon \,, 
$$
so that, if (\ref{persistent}) holds, the probability of contradiction is controlled as evidence accumulates. 

Confidence sequences satisfying (\ref{persistent}) are obtained in Robbins (1970, see formula (3)). A heuristic argument for their consistency is outlined  in the Appendix.  Robbins' regions,    denoted by $\hat{\Theta}_{1-\varepsilon}( Y^{(n)})$, with realization $\hat{\Theta}_{1-\varepsilon}( y^{(n)})$,   have the form 
\begin{equation}\label{conf_reg}
\hat{\Theta}_{1-\varepsilon}( y^{(n)}) = \left\{ \theta \in \Theta \, : \; p_n(y^{(n)};\theta ) \geq\varepsilon  q_n(y^{(n)})  \right\}\,, 
\end{equation}
where $q_n(y^{(n)}) $  is the {mixture  density} 
\begin{equation}\label{mix}
q_n(y^{(n)}) = \int_\Theta p_n(y^{(n)};\theta ) \pi(\theta) \, d\theta \,. 
\end{equation}
In (\ref{mix}), the weight function $\pi(\theta)$ is a preset probability density over $\Theta$ with $\pi(\theta) >0$ for every $\theta \in \Theta$. 
 The value $1-\varepsilon$ will be called here the  persistence level of the confidence sequence (\ref{conf_reg}).   

To illustrate the simplicity  of the approach,  
the proof in Robbins (1970) that the sequence of regions $\hat{\Theta}_{1-\varepsilon}( Y^{(n)})$ satisfies (\ref{persistent}) is sketched in the Appendix. The key argument is inequality  (\ref{eq1}),  giving a bound on the probability of reaching strongly misleading evidence from the likelihood ratio statistic (Royall, 1997, page 7). {The proof does not require the components of $Y^{(n)}$ to be independent or identically distributed.}


While conventional inference --- both Bayesian and frequentist ---  is contingent on the current sample or on the  generating mechanism of the current sample,  inference from Robbins' confidence sequences is in a sense enduring; indeed, it leads to conclusions that with reasonably high probability can withstand any further scrutiny under the same data generating model.

Robbins' confidence sequences are likelihood-based.  Specifically, $\hat{\Theta}_{1-\varepsilon}( y^{(n)})$ is the region  of $\theta$ values whose likelihood 
$L(\theta; y^{(n)})=p_n(y^{(n)};\theta )$ is  larger than a fixed fraction of the {mixture density}  $ q_n(y^{(n)})$.  Therefore, regions  (\ref{conf_reg}) are invariant under one-to-one transformations of $y$ and one-to-one transformations of $\theta$. 
The mixture density  $ q_n(y^{(n)})$ incorporates prior information, possibly notional.   In any case,   the importance of the specification  of $\pi(\theta)$ is 
downplayed    because property  (\ref{persistent}) holds for every  strictly positive $\pi(\theta)$.
{Confidence regions (\ref{conf_reg}) are nested, i.e., $\hat{\Theta}_{1-\varepsilon'}( y^{(n)})\subseteq \hat{\Theta}_{1-\varepsilon}( y^{(n)}) $, when  $1-\varepsilon'< 1-\varepsilon$. The maximum likelihood estimate $\hat\theta_n$ is always in $\hat{\Theta}_{1-\varepsilon}( y^{(n)})$, being
$p_n(y^{(n)};\hat\theta_n ) \geq  q_n(y^{(n)})$. As a first illustration, in Example 1 below we obtain  Robbins' confidence sequence for the mean of a normal distribution. Further examples and simulation results are given in Section 4.}\\

{
{\it Example 1. {Inference about the mean of a normal population, known variance.}}\\
Suppose that $Y_i, \, i=1,2, \ldots $, are i.i.d.\ $N(\theta , \sigma^2_0)$, with unknown mean $\theta$ and known variance $\sigma^2_0$. 
Reduction by sufficiency produces the sequence  of sample means $\bar{Y}_n = \sum_{i=1}^n Y_i/n$ with model $N(\theta , \sigma^2_0/n)$, $n=1,2, \ldots$. 
The density of $\bar{Y}_n$ under  $\theta$ is  
$$
p_n(\bar{y}_n ; \theta) = \dfrac{\sqrt{n}}{\sqrt{2\pi \sigma_0^2}} \exp \left\{- \dfrac{n ( \bar{y}_n - \theta )^2}{2\sigma_0^2} \right\} \,. 
$$ 
Taking  as a weight function the $N(\mu_0 , \tau_0^2)$ density,  
$
\pi (\theta) = (2\pi \tau_0^2)^{-1/2} \exp \left\{ -(2\tau_0^2)^{-1}(\theta - \mu_0)^2  \right\}$,
the mixture distribution of $\bar{Y}_n$ is  $N(\mu_0 , \tau_0^2 + \sigma^2_0/n)$, giving  
$$
q_n(\bar{y}_n ) = \dfrac{1}{\sqrt{2\pi} \sqrt{ \tau^2_0 + \sigma_0^2/n}} \exp \left\{- \dfrac{1}{2} \dfrac{( \bar{y}_n - \mu_0 )^2}{\tau^2_0 +\sigma_0^2/n }\right\} \,. 
$$
After some algebra, Robbins' confidence sequence 
$
\hat{\Theta}_{1-\varepsilon}(\bar{y}_n) = \left\{ \theta \in \reali \, : \;  p_n(\bar{y}_n ; \theta) \geq \varepsilon  q_n(\bar{y}_n ) \right\}
$
is seen to consist of the intervals   $\bar{y}_n \pm  d_n(\sigma^2_0)$, where
\begin{equation}\label{dn}
d_n(\sigma^2_0)= \dfrac{\sigma_0}{\sqrt{n}} \sqrt{\log \dfrac{\tau_0^2+\sigma^2_0/n}{\sigma_0^2/n} + \dfrac{(\bar{y}_n-\mu_0)^2}{\tau_0^2+\sigma^2_0/n}  -2\log \varepsilon  } \, . 
\end{equation}
Numerical evaluation of contradictions and non-coverages of the confidence sequence $\bar{y}_n \pm  d_n(\sigma^2_0)$ is given in Example 2. }\\

When the parameter is partitioned as $\theta = (\psi , \lambda)$, where $\psi \in \Psi$ is a $p_0$-dimensional component of interest and $\lambda $ is nuisance, in some cases inference about $\psi$ can be based on a statistic $t^{(n)}=t(y^{n})$ producing a  marginal or conditional model  free of $\lambda$. In these cases,  $p_n (t^{(n)}; \psi)$ or $p_n(y^{(n)} | t^{(n)} ; \psi)$ may replace $p_n (y^{(n)} ; \theta)$ in (\ref{conf_reg}) with $q(y^{(n)})$  redefined accordingly.  However,  Robbins' confidence sequences for a parameter of interest are also obtainable when a reduction by marginalization or conditioning is not available, without requiring the calculation of sampling distributions.  {Indeed, a confidence sequence for $\psi$ with persistence level $1-\varepsilon$ is given by} the projection of $\hat{\Theta}_{1-\varepsilon} (y^{(n)})$ on $\Psi$, 
$$
\hat{\Psi}_{1-\varepsilon}( y^{(n)}) = \left\{ \psi \in \Psi \, : \;  (\psi , \lambda) \in \hat{\Theta}_{1-\varepsilon}( y^{(n)}) \text{ for some } \lambda  \right\}
$$
and can be expressed in terms of the profile likelihood $p_n(y^{(n)};\psi, \hat{\lambda}_\psi ) $ as 
\begin{equation}\label{profconfseq}
\hat{\Psi}_{1-\varepsilon}( y^{(n)}) = \left\{ \psi \in \Psi \, : \;  p_n(y^{(n)};\psi, \hat{\lambda}_\psi ) \geq \varepsilon  q_n(y^{(n)})  \right\}\,.
\end{equation}
Above, $\hat{\lambda}_\psi$ is the maximum likelihood estimate of $\lambda$ in the model for $y^{(n)}$ with $\psi$ fixed 
and 
$
q_n(y^{(n)}) 
$ is given by (\ref{mix}).

{Sequences (\ref{profconfseq}) are likely to be far more conservative than their counterpart with known $\lambda$, as discussed in Example 2.  Precise quantification of the nuisance parameters effect --- for instance in the spirit of DiCiccio et al.\ (2015) --- seems to be out of reach. }

{On the practical side, suppose that a normal approximation is available for the maximum likelihood estimator $\hat\psi_n$ of a scalar $\psi$, i.e.\ $\hat\psi_n \stackrel{\cdot}{\sim} N(\psi, v_n)$, with $v_n$ an estimate of the asymptotic variance of  $\hat\psi_n$. If   
a $N(\psi_0 , \tau_0^2)$ density is used  as a weight function,  a closed form approximate confidence sequence for $\psi$ is  
\begin{equation}\label{hatdn} 
\hat\psi_n \pm d_n(n v_n)\,,
\end{equation}
 where  $d_n(\cdot)$ is given by (\ref{dn}).  It will be seen through simulations in Section 4   that this proposal seems to maintain the persistence level $1-\varepsilon$ in all the examples considered. Intervals (\ref{hatdn})  have a Wald-type structure, so that they are no longer  exactly equivariant under reparameterizations.} \\

\section{Frequentist, pure likelihood and Bayesian properties}
When $n$ is large, the coverage probability of $\hat{\Theta}_{1-\varepsilon}( Y^{(n)})$ is close to one. In fact, from the usual asymptotics where, under $\theta$,  
$$
2\left( \ell(\hat\theta_n;Y^{(n)} ) - \ell(\theta;Y^{(n)} ) \right) \stackrel{d}{\to} \chi^2_p\,,
$$
we have 
\begin{equation}\label{res}
\lim_{n \to \infty} P_\theta ( \theta \in \hat{\Theta}_{1-\varepsilon}( Y^{(n)})) =1 \,. 
\end{equation}
Details are given in the Appendix. {Thus, the confidence level of $\hat{\Theta}_{1-\varepsilon}( Y^{(n)})$ is adjusted to the sample size}. 
This behaviour contrasts greatly with what is usually sought for in frequentist inference, i.e., asymptotic coverage equal to the nominal level, $1 -\alpha$.  
 Under this respect, a frequentist statistician willing to ensure her confidence regions to be non-contradictory with  positive probability  has to pay a  price in terms of overcoverage for fixed $n$. {For fixed $n$, this of course implies, a larger probability of covering a given false parameter value.}  
 
%

Robbins' confidence sequences entail a novel concept of confidence, involving the current size $n$ experiment and its future enlargements, hypothetical or not. A persistence level $1-\varepsilon$ ensures that
\begin{equation}\label{ensure} 
P_\theta\left( \theta \in \hat{\Theta}_{1-\varepsilon}( Y^{(m)})  \;\;{\rm for\;\; every}\;\; m \geq n  \right) \geq 1 -\varepsilon \, , 
\end{equation}
so that  
$$
P_\theta\left( \cap_{m \geq n} \hat{\Theta}_{1-\varepsilon}( Y^{(m)})  \neq \emptyset  \right) \geq 1 -\varepsilon \, .  
$$
In practice, we have high confidence that no contradiction with the current conclusions  $\hat{\Theta}_{1-\varepsilon}( y^{(n)})$ would occur with larger sample sizes, even in settings where the sample enlargement is only hypothetical.

It is important to stress that what happened for sample sizes  from 1 to $n-1$ does not matter. Moreover, although  the sequence $ \cap_{j \leq n} \hat{\Theta}_{1-\varepsilon}( Y^{(j)})$ satisfies (\ref{persistent}) as well (see Robbins, 1970, Section 3), it is not eligible as a  sensible confidence sequence because  $\cap_{j \leq n} \hat{\Theta}_{1-\varepsilon}( y^{(j)})$ could be empty,  and therefore not consistent.  


From 
\begin{eqnarray*}
P_\theta \left( \theta \in \cap_{m \geq n} \hat{\Theta}_{1-\varepsilon}( Y^{(m)}) \right)  &=&  P_\theta \left( \theta \in \hat{\Theta}_{1-\varepsilon}( Y^{(n)}) \right)\\
&&    \qquad P_\theta \left( \theta \in \cap_{m > n} \hat{\Theta}_{1-\varepsilon}( Y^{(m)}) \mid \theta \in \hat{\Theta}_{1-\varepsilon}( Y^{(n)}) \right)\,,  
\end{eqnarray*}
the frequentist assurance of $\hat{\Theta}_{1-\varepsilon}( y^{(n)})$ expressed by (\ref{ensure}) entails  
$$
 P_\theta \left( \theta \in \cap_{m > n} \hat{\Theta}_{1-\varepsilon}( Y^{(m)}) \mid \theta \in \hat{\Theta}_{1-\varepsilon}( Y^{(n)}) \right)  
= \dfrac{P_\theta \left( \theta \in \cap_{m \geq n} \hat{\Theta}_{1-\varepsilon}( Y^{(m)} \right)}{P_\theta \left( \theta \in \hat{\Theta}_{1-\varepsilon}( Y^{(n)}) \right)}\\
 \geq 1-\varepsilon \,.
$$
Therefore, if $\hat{\Theta}_{1-\varepsilon}( y^{(n)})$ covers the truth --- an easily conceded premise if $n$ is large enough, in view of (\ref{res}) --- then, with 
probability at least $1-\varepsilon$,  no contradiction with $\hat{\Theta}_{1-\varepsilon}( y^{(n)})$ will be seen under future enlargements of the study.

Moreover, the reward for  overcoverage in a fixed $n$ perspective is that 
Robbins' confidence sequence (\ref{conf_reg}) offers inference that rarely fails to reproduce even in a multiple investigation perspective. 
Let the sequences $Y^{(n)}$ and $Y^{*(n')}$ be independent with 
the same statistical model $\left\{ P_\theta, \; \theta \in \Theta \subseteq \reali^p \right\}$ and the same true parameter value.  
Statistician A will observe the initial part of the  sequence $Y^{(n)}$, statistician B will observe the initial part of the sequence $Y^{*(n')}$. 
If both A and B adopt and communicate publicly Robbins' confidence regions with the same $\varepsilon$, though with possibly different preset weight functions, they will be usually found in 
agreement,  because   
\begin{eqnarray*}
&&P_\theta \left(  \hat{\Theta}_{1-\varepsilon}( Y^{(n)}) \cap \hat{\Theta}_{1-\varepsilon}( Y^{*(n')}) \neq \emptyset  \right) \\&&\geq
P_\theta \left( \theta \in \hat{\Theta}_{1-\varepsilon}( Y^{(n)}) \cap \hat{\Theta}_{1-\varepsilon}( Y^{*(n')}) \quad {\rm for\;\; every}\;\; n, n' \geq 1 \right) \geq (1-\varepsilon)^2 \,. 
\end{eqnarray*}

Regions $\hat{\Theta}_{1-\varepsilon}( y^{(n)}) $, depending on the data only through the likelihood function,  agree with the strong likelihood principle. As a consequence, they are insensitive to the stopping rule and can be used when the stopping rule is unknown. Moreover, Robbins' confidence sequences obey  both the sufficiency and the conditionality principles. 
For sufficiency, let $s^{(n)}=s(y^{(n)})$ be a sufficient statistic for $p_n(y^{(n)};\theta)$, $\theta\in\Theta\subseteq \reali^p$, so that $$p_n(y^{(n)};\theta)=p_{S^{(n)}}(s^{(n)};\theta)\,p_n(y^{(n)}|s^{(n)}),$$ with $p_{S^{(n)}}(s^{(n)};\theta)$ the marginal density of   $S^{(n)}= s(Y^{(n)})$ and    $p_n(y^{(n)}|s^{(n)})$ the conditional density of $Y^{(n)}$ given $S^{(n)}=s^{(n)}$.    Then,  $\hat{\Theta}_{1-\varepsilon}( y^{(n)})=
\hat{\Theta}_{1-\varepsilon}( s^{(n)})$. 
As to conditionality, let $a^{(n)}=a(y^{(n)})$ be a distribution constant statistic, 
so that $$p_n(y^{(n)};\theta)=p_{A^{(n)}}(a^{(n)})\,p_n(y^{(n)}|a^{(n)};\theta).$$
Then 
$$\hat{\Theta}_{1-\varepsilon}( y^{(n)})=
\left\{\theta\in\Theta\,:\;    p_n(y^{(n)}|a^{(n)};\theta) 
\geq
\varepsilon \int_{\Theta} p_n(y^{(n)}|a^{(n)};\theta)  \pi(\theta) d\theta
\right\}  = \hat{\Theta}_{1-\varepsilon}( y^{(n)}|a^{(n)})\,,
$$
so that regions $\hat{\Theta}_{1-\varepsilon}( Y^{(n)})$ have probability of contradiction bounded by $\varepsilon$ also conditionally on $a^{(n)}$. 


Let us consider now Bayesian properties of the confidence sequence (\ref{conf_reg}). 
%
When  the mixing density $\pi(\theta)$ represents a prior distribution,   
 the posterior with data $y^{(n)}$  is 
$$
\pi(\theta|y^{(n)}) = \frac{p_n(y^{(n)};\theta ) \pi(\theta)}{ \int_\Theta p_n(y^{(n)};\theta ) \pi(\theta) \, d\theta  } \,.
$$ 
Definition (\ref{conf_reg}) may be recast as
\begin{equation}\label{repr_Bayes}
\hat{\Theta}_{1-\varepsilon}( y^{(n)}) = \left\{ \theta \in \Theta \, : \; \pi(\theta|y^{(n)}) \geq \varepsilon  \pi(\theta)  \right\}\,. 
\end{equation}
The complementary set $\bar{\Theta}_{1-\varepsilon}(y^{(n)})=\Theta\setminus \hat{\Theta}_{1-\varepsilon}( y^{(n)})$ has thus posterior probability 
$$
  \int_{\bar{\Theta}_{1-\varepsilon}(y^{(n)})} \pi(\theta|y^{(n)})\, d\theta \leq \varepsilon \int_{\bar{\Theta}_{1-\varepsilon}(y^{(n)})} \pi(\theta)\, d\theta \leq \varepsilon \,,
$$
so that $\hat{\Theta}_{1-\varepsilon}( y^{(n)})$ has posterior probability greater than  $1-\varepsilon$. 
{Representation (\ref{repr_Bayes}) shows that inference from Robbins' confidence sequences proceeds by subtraction as the posterior concentrates around the true parameter value, eliminating from $\Theta$ the most implausible values. }

Unlike the usual credible regions $\hat{\Theta}_{B}^{1-\alpha}( y^{(n)})$ satisfying
$$
\int_{\hat{\Theta}_{B}^{1-\alpha}( y^{(n)})}
\pi(\theta|y^{(n)}) \, d\theta =1-\alpha\,,
$$
whose credibility $1-\alpha$ is not adjusted to the sample size, regions $\hat{\Theta}_{1-\varepsilon}( y^{(n)})$ have bounded probability of being contradictory even in a Bayesian sense. Indeed, let $P$ be the joint probability model of $\theta$ and $Y^{(\infty)}$, where $\theta$ has marginal density $\pi(\theta)$ and, given $\theta$,  $Y^{(\infty)}$ has conditional distribution $P_{\theta}$.  
In this setting, (\ref{persistent})  is a conditional probability statement. With $\hat{\Theta}_n=\hat{\Theta}_{1-\varepsilon}( Y^{(n)})$,  formula (\ref{persistent}) implies that
$$
P\left( \theta \in \bigcap_{n\geq 1}  \hat{\Theta}_{1-\varepsilon}( Y^{(n)})
 \right) \geq 1 -\varepsilon \, .
$$

%
  
 \section{Examples}

The implementation of Robbins' confidence sequences requires the specification of $\pi(\theta)$ and the choice of $\varepsilon$ or a range of $\varepsilon$ values. {These issues, together with an assessment of the approximate confidence sequences (\ref{hatdn}), are illustrated through the following examples. In each scenario, empirical percentages of contradictions and non-coverages are evaluated through simulation over a range $n_{min} \leq n \leq n_{max}$ of sample sizes. Let $ ( \underline{\theta}_n  , \bar{\theta}_n)$ be a sequence of confidence intervals  for a scalar parameter $\theta$ and denote by  $\min$ and $\max$ the minimum and maximum over the range of interest.  A sequence shows a contradiction whenever $\max (\underline{\theta}_n) > \min (\bar{\theta}_n)$ and a  non-coverage of $\theta$  whenever $\max (\underline{\theta}_n) >\theta$ or $\min (\bar{\theta}_n) <\theta$.    
}\\

{\it Example 2. {Inference about the mean of a normal population.}}\\
In the setting of Example 1, a  simulation study has been performed in order to compare properties of  Robbins' confidence sequences with those of customary confidence intervals.   
Contradictions and non-coverages have been monitored for 10,000 replications of enlarging samples of size $n$ with $n_{min}=10$ and $n_{max}= 4,000$.  

In Table 1 the simulation results for confidence intervals $\bar{y}_n \pm \sigma_0 z_{1-\alpha/2}/\sqrt{n}$, with  confidence level  $1-\alpha = 0.90, 0.95,0.99,0.995$, are shown. 
Contradictions and non-coverages are dominant for the 90\%  and 95\%  levels.  They are both comparatively uncommon for the level 99.5\%,  but their relative frequency could be made as close to 1 as desired by letting $n_{max}$ large enough, when all non-coverages become contradictions. The simulation has been performed by sampling  standard normal deviates. The results, however, do not depend on the true value of the parameters of the normal population.

\begin{table}[h!]
\begin{center}
\caption{Normal population with known variance 1: empirical percentages of contradictions and non-coverages for some $n$ of  intervals for the mean with confidence level $1-\alpha$ in 10,000 sequences of samples   with size   from 10 to 4,000. The true value of $\theta$ is 0.}\label{Table 1}
\begin{tabular}{cccccc}
&&&&&\\
$100(1-\alpha)$ & &   90    &   95    &   99     &     99.5  \\
\hline
&&&&&\\
   contradictions  & & 51.32  &  27.35  &  5.20   &    2.29     \\
   non-coverages   &  &   77.79  &  54.21  &  18.52  &   10.86   \\ 
&&&&&\\
\hline
        \end{tabular}
\end{center}
\end{table}

Table 2 displays the  simulation  results for Robbins' confidence sequences with persistence levels $1-\varepsilon = 0.50, 0.80, 0.90, 0.95$ and various  $N(\mu_0 , \tau_0^2)$ densities as a weight function $\pi(\theta)$ for $\theta$. 
When $\pi(\theta)$ is concentrated around the true $\theta$, contradictions and non-coverages are comparatively abundant, but their relative frequency remains under the threshold $\varepsilon$.  When  $\pi(\theta)$ is discrepant from the likelihood, that is  $\mu_0$ is far from the true $\theta$, 
 the conflict between the weight function and the likelihood is resolved in favour of the likelihood, through wider confidence intervals.  
This counterbalance increases conservativeness of the sequence with respect to the $\varepsilon$ bound. Apart from these cases,  when $1-\varepsilon =0.80$ the results in terms of observed contradictions and non-coverages for some $n$ in the range 10--4,000 are qualitatively comparable with those for the customary intervals with confidence level 0.995. {As a numerical illustration, the confidence interval with $1-\alpha=0.995$ when $n=100$, $\sigma^2_0=1$ and $\bar{y}_{100}=0$, is $\pm 0.281$, while Robbins' confidence intervals with $1-\varepsilon=0.80$ are 
$\pm 0.280$ when $(\mu_0,\tau^2_0)=(0,1)$,  
$\pm 0.304$ when $(\mu_0,\tau^2_0)=(0,4)$, 
$\pm 0.297$ when $(\mu_0,\tau^2_0)=(1,1)$, 
$\pm 0.308$ when $(\mu_0,\tau^2_0)=(1,4)$.}

{If also the variance is unknown,  the model has parameter $\theta=(\mu,\sigma^2)$ and a confidence sequence for $\mu$ may be obtained from (\ref{profconfseq}). A convenient form of the weight function $\pi(\theta)$ is that of a normal-inverse gamma  conjugate prior, where $1/\sigma^2$ has a gamma distribution with shape parameter $\alpha_0$ and rate $\beta_0$, and, conditionally on $\sigma^2$, $\mu$ has a normal distribution with mean $\mu_0$ and variance $\sigma^2/\kappa_0$. With this specification, $q_n(y^{(n)})$ has a closed form expression and the confidence sequence   (\ref{profconfseq}) is of the form $\bar{y}_{n}\pm \hat\sigma_n h_n$, where $\hat\sigma^2_n$ is the maximum likelihood estimate of $\sigma^2$ and $h_n$ is an explicit function of $n$, $\alpha_0$, $\beta_0$, $\mu_0$, $\kappa_0$, $\bar{y}_n$, $\hat\sigma^2_n$. }

{When $n$ is  sufficiently large, intervals (\ref{hatdn}) with $\hat\psi=\bar{y}_{n}$ and $v_n= \hat\sigma^2_n/n$ may be considered as a simple approximate solution. Simulation results, not reported here, with  $n_{min}=30$ and $n_{max}= 4,000$, indicate that confidence sequences (\ref{profconfseq}) are much more conservative than confidence sequences (\ref{hatdn}) computed with the normal weight functions having the same mean and variance as the marginal conjugate for $\mu$. Moreover, with the same range of sample sizes, intervals (\ref{hatdn})  show empirical percentages of contradictions and non-coverages only slightly larger than their counterparts with known $\sigma^2$. Continuing the previous numerical illustration, and assuming $\hat\sigma^2_n=1$, $\mu_0=1$, $\kappa_0=8$, $\alpha_0=2$, $\beta_0=1$,  we get, for interval $\bar{y}_{n}\pm \hat\sigma_n h_n$ from (\ref{profconfseq}),  $h_n=0.323$, while, for interval (\ref{hatdn}), using the corresponding  $N(0,0.125)$ weight function, $d_n(\hat\sigma^2_n)=0.241$. 
}  \\

\begin{table}[h!]
\begin{center}
\caption{Normal population with known variance 1: empirical percentages of contradictions and non-coverages   for Robbins' confidence sequences for the mean with persistence level $1-\varepsilon$ in 10,000 sequences of samples with size  from 10 to 4,000  and various normal weight functions.  The true value of $\theta$ is 0.}\label{Table 2}
\begin{tabular}{ccccccc}
&&&&&&\\
    $N(\mu_0 , \tau_0^2)$ weight function  &$100(1-\varepsilon)$     &       &   50     &    80     &   90     &     95    \\
\hline
&&&&&&\\
$\mu_0=0$, $\tau_0^2 =0.1$ &      contradictions      &       &  17.62  &    4.08  &    1.31   &     0.48     \\
                                                &      non-coverages       &       &  39.06   &   15.31  &    7.28   &     3.81     \\ 
&&&&&&\\
$\mu_0=0$, $\tau_0^2 =1.0$ &      contradictions      &       &  10.35   &    3.21  &    1.37   &    0.59     \\
                                                &      non-coverages       &       &  22.05   &   9.39  &    4.68   &     2.30     \\ 
&&&&&&\\
$\mu_0=0$, $\tau_0^2 =10$  &      contradictions      &       &  3.03   &    1.06  &    0.48   &   0.23     \\
                                                &      non-coverages       &       &  8.42   &     3.38  &    1.69   &    0.97     \\ 
&&&&&&\\
$\mu_0=1$, $\tau_0^2 =1.0$ &      contradictions      &       &  6.54   &     2.03  &    0.96   &  0.46     \\
                                                &      non-coverages       &       &  14.75   &    6.15  &    3.02   &  1.59     \\ 
&&&&&&\\
$\mu_0=2$, $\tau_0^2 =1.0$ &      contradictions      &       &  1.89   &    0.72  &    0.34   &     0.20     \\
                                                &      non-coverages       &       &  4.75   &     1.97  &    1.09  &     0.71     \\ 
&&&&&&\\
$\mu_0=5$, $\tau_0^2 =1.0$ &      contradictions      &       &    0.00  &     0.00  &    0.00   &     0.00     \\
                                                &      non-coverages       &       &    0.01   &     0.00  &    0.00   &     0.00     \\ 
&&&&&&\\
\hline
\end{tabular}
\end{center}
\end{table}

\noindent
{\it Example 3: Bernoulli population.}\\
Suppose that $Y_i, \, i=1,2, \ldots $ are i.i.d.\ Bernoulli $Bi(1, \theta )$, with unknown mean $\theta \in (0,1)$. 
Reduction by sufficiency produces the sequence of sample sums $S_n = \sum_{i=1}^n Y_i$, whose model is $Bi(n, \theta )$, with density under $\theta$   
$$
p_n(s_n ; \theta) = {n \choose s_n} \theta^{s_n} (1-\theta )^{n - s_n} \,. 
$$ 
Let us consider as a weight function a conjugate $Beta(\alpha , \beta)$ density 
$$
\pi (\theta) = \dfrac{1}{B(\alpha , \beta)} \theta^{\alpha -1} (1-\theta)^{\beta -1}\,,  
$$
where $\alpha, \beta >0$ and $B(\alpha , \beta)= \Gamma (\alpha)\Gamma (\beta )/\Gamma(\alpha + \beta)$.  
The mixture distribution of $S_n$ is  then beta-binomial,  with density  
$$
q_n(s_n ) =  {n \choose s_n} \dfrac{B(s_n+\alpha, n-s_n + \beta)}{B(\alpha , \beta)} \,. 
$$
The choice $\alpha=\beta=0.5$ corresponds to Jeffreys' prior. When  $\alpha=\beta =1$ the weight function is a continuous uniform distribution on $[0,1]$ and the  mixture distribution of $S_n$ is   discrete uniform on $\{ 0, 1, \ldots , n  \}$. 
Robbins' confidence sequence  
$$
\hat{\Theta}_{1-\varepsilon}(s_n) = \left\{ \theta \in \reali \, : \;  p_n(s_n ; \theta) \geq \varepsilon  q_n(s_n) \right\}
$$
does not have a closed-form expression but can be easily computed numerically{, because the log likelihood function is concave}.  

Intervals with asymptotic confidence level $1-\alpha$ from the likelihood ratio statistic  have the form 
$$
\tilde{\Theta}_{1-\alpha}(s_n)= \left\{ \theta \in (0,1) \, : \; p_n(s_n; \theta) \geq p_n (s_n; \hat\theta_n) \exp\{-0.5 \chi^2_{1,1-\alpha} \} \right\}\,,
$$   
where $\hat\theta_n = \bar{y}_n=s_n/n$ is the maximum likelihood estimate and $\chi^2_{1,1-\alpha}$ is the $(1-\alpha)$-quantile of a chi-squared distribution with 1 degree of freedom.  {Also} these intervals are easily computed numerically.  

A small simulation study with various true $\theta$ values and various weight functions has been performed. In particular, contradictions and non-coverages for $n$ in the range with $n_{min}=100$ and $n_{max}= 4,000$ have been enquired. {Here we considered  $n_{min}=100$ to rely on standard asymptotics of the likelihood ratio statistic.}  The number of replications remains $10,000$.  

Table \ref{Table 3} displays the simulation results for the confidence intervals with confidence level $1-\alpha = 0.90, 0.95, 0.99, 0.995$ obtained from the likelihood ratio statistic. Contradictions and non-coverages are important when $1-\alpha=0.90, 0.95$. The level $1-\alpha=0.995$ gives  a marked improvement.

\begin{table}[h!]
\begin{center}
\caption{Bernoulli population: empirical percentages of contradictions and non-coverages for some $n$ of likelihood  ratio  intervals for the mean with confidence level $1-\alpha$ in 10,000 sequences of samples with size  from 100 to 4,000.}\label{Table 3}
\begin{tabular}{ccccccc}
&&&&&&\\
                         & $100(1-\alpha)$ &   &   90    &   95     &   99     &     99.5  \\
\hline
&&&&&&\\
$\theta = 0.5$  &   contradictions   &   &   28.38   &  12.06    &   1.42     &   0.57     \\  
                         &   non-coverages    &   &   64.44  &   41.78  &   12.58     &    7.23     \\  
&&&&&&\\
$\theta = 0.7$  &   contradictions   &   &   27.83      &   12.07     &   1.54   &     0.54     \\ 
                         &   non-coverages    &   &   64.47   &  42.94    &   13.15   &     7.58     \\ 
&&&&&&\\
$\theta = 0.9$  &   contradictions   &   &   28.24     &   12.16     &   1.60      &     0.65   \\ 
                         &   non-coverages    &   &   62.79     &   41.40   &   12.66      &     7.31   \\ 
&&&&&&\\
\hline
\end{tabular}
\end{center}
\end{table}

\begin{table}[h!]
\begin{center}
\caption{Bernoulli population: empirical percentages of contradictions and non-coverages for some $n$ of Robbins' confidence sequences for the mean with  persistence level $1-\varepsilon$ in 10,000 sequences of samples with size from 100 to 4,000   and various  weight functions on the mean.}\label{Table 4}
\begin{small}
\begin{tabular}{cccccccc}
&&&&&&&\\
 true $\theta$          &    weight function                                 &$100(1-\varepsilon)$     &       &   50     &    80     &   90     &     95    \\
\hline
&&&&&&&\\
          0.5                 & $Beta(.5,.5)$     &      contradictions      &       &  0.86   &     0.17  &    0.04   &     0.01     \\   
                                &                          &      non-coverages       &       &  7.36   &    3.25  &    1.47   &     0.75     \\    
&&&&&&&\\
          0.5                 &  $Beta(1,1)$      &      contradictions      &       &  1.54   &     0.36  &    0.10   &     0.03     \\    
                                &                          &      non-coverages       &       &  10.85   &    4.73  &    2.46   &     1.12     \\     
&&&&&&&\\
          0.5                 &  $Beta(5,5)$      &      contradictions      &       &  4.82   &     1.07    &    0.33    &     0.08     \\   
                                &                          &      non-coverages       &       &  21.42   &    9.29    &    4.97   &     2.46     \\  
&&&&&&&\\
          0.7                 & $Beta(.5,.5)$     &      contradictions      &       &  0.84    &     0.26   &    0.08   &    0.02     \\ 
                                &                          &      non-coverages       &       &  7.29   &     3.27   &    1.70   &    0.88     \\   
&&&&&&&\\
          0.7                 &  $Beta(1,1)$      &      contradictions      &       &  1.40   &     0.40  &    0.11   &     0.03     \\  
                                &                          &      non-coverages       &       &  9.87   &    4.28  &    2.33   &     1.31     \\  
&&&&&&&\\
          0.7                 &  $Beta(5,5)$      &      contradictions      &       &  2.21   &    0.54  &    0.18   &    0.05     \\ 
                                &                          &      non-coverages       &       &  11.63   &    5.26  &    2.75   &     1.49     \\  
&&&&&&&\\
 0.9                 &  $Beta(.5,.5)$      &      contradictions      &       &  0.86   &     0.22  &    0.08   &     0.04     \\  
                                &                          &      non-coverages       &       &  7.05   &    2.94 &    1.49   &    0.69     \\  
&&&&&&&\\
          0.9                 &  $Beta(1,1)$      &      contradictions      &       &  0.69   &     0.14  &    0.04   &     0.02      \\ 
                                &                          &      non-coverages       &       &  6.14   &    2.60   &    1.27   &     0.59    \\ 
&&&&&&&\\
                   0.9                 & $Beta(5,5)$     &      contradictions      &       &  0.06   &     0.02  &    0.01   &     0.01     \\ 
                                &                          &      non-coverages       &       &  1.00   &    0.35  &    0.22   &     0.13     \\   

&&&&&&&\\
\hline
\end{tabular}
\end{small}
\end{center}
\end{table}

In Table \ref{Table 4}  results for Robbins' confidence sequences with persistence levels $1-\varepsilon = $ 0.50, 0.80, 0.90, 0.95 
and various beta weight functions are shown. 
When the weight function is centered at the true $\theta$,  non-coverages are comparatively abundant.
Contradictions are rarely observed in the range of $n$ values considered. As expected,  conservativeness increases as the weight function moves away from the true parameter value. Again, the results 
 for $
\hat{\Theta}_{1-\varepsilon}(s_n) $ when $1-\varepsilon =0.80$ are qualitatively comparable with those for $
\tilde{\Theta}_{1-\alpha}(s_n)$  with $1-\alpha=0.995$. 
 
{When $n_{min}$ is sufficiently large, an approximate confidence sequence for $\theta$ that does not require numerical calculation is obtained using (\ref{hatdn}) in the variance stabilizing parameterization $\omega(\theta)=\arcsin\sqrt \theta$, and with a $N(\mu_0, \tau^2_0)$ weight function for $\omega$. Since $\omega(\bar{Y}_n)$ is approximately distributed as $N(\omega(\theta), 1/(4n))$,  an approximate confidence sequence for $\omega(\theta)$ is 
\begin{equation}\label{bin_arcsin}
\omega(\bar{y}_n)\pm \dfrac{1}{2\sqrt{n}} \sqrt{ \log \dfrac{\tau_0^2+1/(4n)}{1/(4n)} + \dfrac{(\omega(\bar{y}_n)-\mu_0)^2}{\tau_0^2+1/(4n)}  -2\log \varepsilon  } \,.
\end{equation}
A simulation study has been done using the same settings as those considered for Table \ref{Table 4}, with   weight functions for $\omega$ chosen as the densities of normal distributions with the same mean and variance as the $\omega(\cdot)$ transformation of the beta weight functions for $\theta$.
Results, not reported here, give empirical percentages of contradictions and non-coverages very similar to those in Table \ref{Table 4}. 
As a numerical illustration, with $n=100$, $\bar{y}_n=0.4$ and the same three beta weights as in Table \ref{Table 4},  Robbins' confidence sequences with $\varepsilon=0.2$, give the intervals
$(0.2673, 0.5435)$, $(0.2738, 0.5359)$ and $(0.2858, 0.5221)$, respectively. The corresponding approximate intervals obtained from (\ref{bin_arcsin}), transformed back in the $\theta$ parameterization, are $(0.2697, 0.5379)$, $(0.2740, 0.5332)$ and $(0.2843, 0.5216)$, giving a quite accurate explicit approximation of Robbins' intervals. For comparison, the likelihood ratio interval with $1-\alpha=0.995$ is $(0.2702, 0.5400)$.} \\

\noindent
{\it Example 4: Two Bernoulli populations.}\\
Suppose that $Y_{1i}$ and $Y_{2i}, \, i=1,2, \ldots, $ are independent Bernoulli $Bi(1, \theta_j)$ with unknown means $\theta_j \in (0,1)$, $j=1,2$.  Consider the log-odds ratio $\psi=\log [\theta_1(1-\theta_2)/\{\theta_2(1-\theta_1)\}]$ as the parameter of interest. With $n_1$ observations from $Bi(1,\theta_1)$ and $n_2$ observations from $Bi(1,\theta_2)$, reduction by sufficiency gives the sample sums $S_{j}=\sum_{i=1}^{n_j} Y_{ji}$, $j=1,2$. A model depending on $\psi$ only is obtained by conditioning on $t^{(n)}$, the observed value of $T^{(n)}=S_{1}+S_{2}$. The conditional density of $S_{1}$ given $t^{(n)}$ is noncentral hypergeometric (McCullagh and Nelder, 1989, Sections 7.3.2 and 7.4.1).
As a weight function for $\psi$, we consider 
\begin{equation}\label{pi_psi}
\pi(\psi)=\psi \exp(\psi/2)/\{\pi^2(\exp(\psi)-1)\}\,,
\end{equation}
when $\psi\neq 0$, and $\pi(0)$ defined by continuity as $ 1/\pi^2$. This is the density of the log-odds ratio when $\theta_1$ and $\theta_2$ have $Beta(0.5, 0.5)$ independent distributions. The numerical calculation of Robbins' confidence sequence for $\psi$ based on the conditional distribution of $S_{1}$ given $t^{(n)}$ may be performed with the aid of the  R package {\tt BiasedUrn} (Fog, 2015).  

Approximate Robbins' confidence sequences  (\ref{hatdn}) using the continuity-corrected quantities
$$
\hat\psi_n=\log\frac{(s_1+0.5)(n_2-s_2+0.5)}{(n_1-s_1+0.5)(s_2+0.5)}
$$
and
$$
v_n=\frac{1}{s_1+0.5}+\frac{1}{n_1-s_1+0.5}+\frac{1}{s_2+0.5}+\frac{1}{n_2-s_2+0.5}
$$
are much simpler to compute and simulation results are given only for them. 
Estimated contradictions and non-coverages   are shown in Table \ref{Table 5}. Persistence levels are $1-\varepsilon =  0.50, 0.80, 0.90, 0.95$  and $n_1=n_2$ range  from 50 to 2,000. The number of replications is  10,000. We set $\theta_1=0.2$ and $\theta_2=0.25$, so that $\psi=-0.288$. Various other pairs $(\theta_1, \theta_2)$ with the same $\psi$ have been considered, leading always to very similar results.  Six pairs $(\mu_0, \tau^2_0)$ for the normal weight function have been used. The pair $(0, 2\pi^2)$ corresponds to the mean and variance of the distribution with density (\ref{pi_psi}) (cf.\ Morris, 1982, Section 4). Empirical percentages of contradictions and non-coverages of the  approximate confidence sequences respect the nominal bounds $100\varepsilon$. On the other hand, the standard asymptotic Wald intervals with nominal confidence level 0.95 show empirical non-coverages  of about 40\%, while empirical non-coverages of Wald intervals with nominal confidence levels 0.99 and 0.995 are about  12\% and 7\%, respectively.  In the range of sample sizes considered, Robbins' confidence sequences with $1-\varepsilon=0.8$ seem intermediate between the standard intervals with $1-\alpha=0.99$ and $1-\alpha=0.995$. 

As a numerical illustration, we compared  conditional and approximate confidence sequences with $n_1=30$, $n_2=70$, $s_1=20$, $s_2=30$. Robbins' confidence interval based on the conditional distribution with $1-\varepsilon = 0.80$ and weight function (\ref{pi_psi}) is $(-0.195, 2.227)$, while the approximate interval is  $(-0.306,  2.211)$ using a $N(0, 2\pi^2)$ weight and $(-0.125,  2.030)$ using a $N(0,1)$ weight. On the other hand, the exact confidence interval with $1-\alpha=0.99, 0.995$, computed using the  R package {\tt exact2x2} (Fay et al., 2018), are $(-0.193,  2.228)$ and $(-0.293,  2.305)$, while Wald confidence intervals with the same levels, $\hat\psi_n\pm z_{1-\alpha/2} \sqrt{v_n}$, are $(-0.204,  2.109)$ and  $(-0.307,  2.213)$. We see that Robbins' approximate interval with  $1-\varepsilon = 0.80$ and $N(0, 2\pi^2)$ weight is in reasonable agreement with the exact interval with $1-\alpha=0.995$.

\begin{table}[h!]
\begin{center}
\caption{Two Bernoulli populations with $\theta_1=0.2$ and $\theta_2=0.25$: empirical percentages of contradictions and non-coverages for approximate Robbins' confidence sequences  with persistence level $1-\varepsilon$ in 10,000 sequences of samples with $n_1=n_2$  from 50 to 2,000  and various normal priors on $\psi$.  The true value of $\psi$ is -0.288}\label{Table 5}
\begin{tabular}{ccccccc}
&&&&&&\\
     weight function                                    &$100(1-\varepsilon)$     &       &   50     &    80     &   90     &     95    \\
\hline
&&&&&&\\
$\mu_0=0$, $\tau_0^2 =2\pi^2$ &      contradictions      &       &  0.80   &   0.22  &    0.02   &    0.01     \\
                                                &      non-coverages       &       &  8.34   &   3.29  &    1.68   &    0.84     \\ 
&&&&&&\\
$\mu_0=0$, $\tau_0^2 =5.0$ &      contradictions      &       &  2.15   &    0.51  &   0.16   &   0.02     \\
                                                &      non-coverages       &       &  15.34   &   6.16  &   3.04   &   1.51    \\ 
&&&&&&\\
$\mu_0=0$, $\tau_0^2 =1.0$  &      contradictions      &       &  4.81   &    0.82  &   0.24   &   0.05     \\
                                                &      non-coverages       &       &  26.06   &   10.75  &  5.35   &  2.61     \\ 
&&&&&&\\
$\mu_0=0$, $\tau_0^2 =0.1$ &      contradictions      &       &  6.60   &    0.59  &   0.05   &  0.01     \\
                                                &      non-coverages       &       &  37.29   &    13.41  &    6.89   &  3.05     \\ 
&&&&&&\\
$\mu_0=1$, $\tau_0^2 =5.0$ &      contradictions      &       &  1.76   &    0.43  &   0.15   &  0.01     \\
                                                &      non-coverages       &       &  13.31   &    5.51  &  2.72   &     1.33     \\ 
&&&&&&\\
$\mu_0=-1$, $\tau_0^2 =5.0$ &      contradictions      &       &   1.96  &     0.44  &    0.13   &     0.02     \\
                                                &      non-coverages       &       &    14.89   &   5.89  &   2.87   &     1.54     \\ 
&&&&&&\\
\hline
\end{tabular}
\end{center}
\end{table}

\section{Concluding remarks}
Herbert E.\ Robbins is mostly acknowledged in Statistics for his path-breaking introduction of empirical Bayes methods, stochastic approximation methods, and contributions to sequential analysis (Lai and Siegmund, 1986), while his proposal of confidence sequences seems to have been largely neglected in the statistical literature. See, however, Gandy and Hahn (2016) where Robbins' confidence sequences provide a tool to keep in check stochastic simulations. Robbins' confidence sequences also inspired repeated confidence intervals (Jennison and Turnbull, 1989), where coverage of the true $\theta$ is required at a finite (typically small)  number of interim analyses of a study. But,  in the discussion of Jennison and Turnbull (1989), Whitehead (1989) points to situations such as long-term epidemiological studies where a fixed number of analyses ``might become a barrier''. This is often the case in modern  applications that routinely deal with large data sets becoming available a little bit at a time due to continuous monitoring. For this reason, some novel attention to Robbins's proposal is currently being paid in the machine learning literature (see, e.g., Johari et al., 2017).  

In this article, we have stressed the link between non-contradiction and coverage along the whole sequence  as the basis for a novel interest in Robbins' confidence sequences.  
These sequences offer durable inferences, satisfying coverage requirements simultaneously for all sample sizes. 
By contrast, inferences stemming from the usual statistical procedures
satisfy coverage requirements separately for any given sample size and  may be called episodic inferences.
The distinction between episodic and sequential environments appears in artificial intelligence, see Russel and Norvig (2010, Section 2.3.2).  
Robbins' confidence sequences strengthen the standards of confidence in a principled way, and, thanks to their frequentist assurance, offer more compelling summarizations of evidence, being also insensitive to the stopping rule.  The price that is paid for controlling for the probability of non-contradiction is that wider  regions are typically needed. This drawback is not, however,  dramatic, as is seen from the numerical illustrations in Section 4.  

We conclude by suggesting some directions for possible extension of the results in the paper.


{Robbins' confidence sequences, based on inequality (\ref{eq1}), require a correctly specified parametric model.  This raises robustness issues. Preliminary simulation results in the same setting as in Example 1, but with data generated from a Student $t$ distribution with moderate degrees of freedom, suggest that persistence may still be under control, provided that $n_{min}$ is large enough, so that the sample mean is approximately normal. As a general strategy, we suggest to use approximations of the form  (\ref{hatdn}) based on asymptotically normal robust estimators of the parameter of interest. }
 
 {The approximate form (\ref{hatdn}) is easily extended to a vector parameter of interest. 
However, the resulting confidence sequence will depend on the parameterization. This could be avoided, using an approximation based on Laplace expansion (see (\ref{Laplace1}) in the Appendix), but sacrificing closed form expressions.}

  {In order to face complex problems, expecially with multidimensional parameters, more work on computational aspects of Robbins' confidence sequences is needed.}

\section*{Appendix}

 
\noindent
{\it Robbins' confidence sequences have the required persistent coverage (Robbins, 1970)}\\
To see that, for regions of the form  (\ref{conf_reg}), inequality (\ref{persistent}) holds
for every $\theta \in \Theta$, consider that  
$$
P_\theta\left( \theta \in \hat{\Theta}_{1-\varepsilon}( Y^{(n)})  \;\;{\rm for\;\; every}\;\; n \geq 1  \right) = 1- P_\theta\left( \theta \not\in \hat{\Theta}_{1-\varepsilon}( Y^{(n)})  \;\;{\rm for\;\; some}\;\; n\geq 1 \right)
$$
and 
$$
P_\theta\left( \theta \not\in \hat{\Theta}_{1-\varepsilon}( Y^{(n)})  \;\;{\rm for\;\; some}\;\; n\geq 1 \right) = P_\theta\left( \frac{q_n(Y^{(n)})}{p_n(Y^{(n)};\theta )} \geq \frac{1}{\varepsilon}\;\;{\rm for\;\; some}\;\; n\geq 1  \right)\, .
$$
The last probability  does not exceed $\varepsilon$ in force of  a fundamental inequality, see (\ref{eq1}) below, for the likelihood ratio statistic. 

{Let $P_\theta$ and $Q$ denote the joint probability distribution of the sequence $Y^{(\infty)}$ corresponding to densities $p_n(y^{(n)}; \theta)$ and $q_n(y^{(n)})$ for $Y^{(n)}$, $n=1,2,\ldots$, respectively.  Then 
\begin{equation}\label{eq1}
P_\theta\left( \frac{q_n(Y^{(n)})}{p_n(Y^{(n)}; \theta)} \geq k \;\;{\rm for\;\; some}\;\; n       \right) \leq \frac{1}{k} \,, 
\end{equation}
for any $k>0$. Robbins' proof of (\ref{eq1}) is as follows.  
Define the stopping time
$$
N=\min \left\{ n \geq 1 \,: \, \frac{q_n(Y^{(n)})}{p_n(Y^{(n)}; \theta)} \geq k \right\} \,,   
$$
when the inequality is satisfied for a finite $n$,  and $N=\infty$ otherwise. 
 Then 
\begin{eqnarray*}
 P_\theta\left( \frac{q_n(Y^{(n)})}{p_n(Y^{(n)}; \theta)} \geq k \;\;{\rm for\;\; some}\;\; n       \right)  &=& P_\theta(N<\infty)\\&=& \sum_{n \geq 1} P_\theta(N=n)
 = \sum_{n \geq 1} \int_{\{y^{(n)}\,:\; N=n\}}p_n(y^{(n)}; \theta)\, dy^{(n)} \\
 & \leq & \sum_{n \geq 1}  \int_{\{y^{(n)}\,:\; N=n\}}\frac{1}{k} q_n(y^{(n)})\, dy^{(n)} \\
 &=& \frac{1}{k} \sum_{n \geq 1} Q(N=n) = \frac{1}{k}Q(N<\infty)\\
 &\leq & \frac{1}{k} \, . 
\end{eqnarray*}
}
Inequality (\ref{eq1}) also follows from a well-known martingale inequality, see e.g.\  Jacod and Protter (2000, Theorem 26.1). \\

\noindent
{\it A heuristic argument for  the consistency of confidence sequences}\\  
 Rigorous proofs of consistency of $\hat{\Theta}_{1-\varepsilon}( y^{(n)})$ when the density of $Y^{(n)}$ belongs to an exponential family are given by Lai (1976)  and Csenki (1979)  for the one-parameter and the multiparameter case, respectively. 
For models whose  likelihood function obeys the usual regularity conditions (see e.g.\ Severini, 2000, Section 3.4), consistency of $\hat{\Theta}_{1-\varepsilon}( y^{(n)})$  may be seen by the following heuristic argument. 

Assume that  $\hat\theta_n$ is the unique maximum of $L(\theta; y^{(n)})$ in an open neighborhood of the true $\theta$. Let  $\ell(\theta; y^{(n)})=\log L(\theta; y^{(n)})$ be the log likelihood function and let 
$j_n(\theta)= j(\theta; y^{(n)})=-\partial^2 \ell(\theta; y^{(n)})/\partial\theta\partial\theta^\top$ be the observed information. Assume moreover that, as under repeated sampling of size $n$, 
$\ell(\theta; Y^{(n)})=O_p(n)$ and $ j(\hat\theta_n; Y^{(n)})$ is positive definite and of order $O_p(n)$.   
Using Laplace expansion, see e.g.\ Barndorff--Nielsen and Cox (1989, Section 3.3), we have
\begin{equation}\label{Laplace}
q_n(y^{(n)}) = \int_\Theta p_n(y^{(n)};\theta ) \pi(\theta) \, d\theta  = p_n(y^{(n)};\hat\theta_n ) \frac{\pi(\hat\theta_n) (2\pi)^{p/2}}{|j_n(\hat\theta_n)|^{1/2}} \{1+ O(n^{-1}) \} \, , 
\end{equation}
so that 
\begin{equation}\label{Laplace1}
\hat{\Theta}_{1-\varepsilon}( y^{(n)}) = \left\{ \theta \in \Theta \, : \; \ell(\theta;y^{(n)} ) >
\ell(\hat\theta_n;y^{(n)} ) +\log\left(  \frac{\varepsilon \pi(\hat\theta_n) (2\pi)^{p/2}}{|j_n(\hat\theta_n)|^{1/2}}\right) +O(n^{-1})
 \right\}\,. 
\end{equation}
Let $k_n=-\log\left(  \varepsilon \pi(\hat\theta_n) (2\pi)^{p/2}/|j_n(\hat\theta_n)|^{1/2}\right)$. Then, for $\theta' \neq \theta$,
$$
  P_\theta ( \theta' \in \hat{\Theta}_{1-\varepsilon}( Y^{(n)})) 
= 
 P_\theta \left(  \ell(\hat\theta_n; Y^{(n)})-\ell(\theta'; Y^{(n)}) <  k_n +O_p(1) \right)\,,
$$
where $ \ell(\hat\theta_n; Y^{(n)})-\ell(\theta'; Y^{(n)})$ is $O_p(n)$ and positive, while 
\begin{equation}\label{kn}
k_n =\frac{p}{2}\log n+O_p(1) \,.
\end{equation}
Therefore, $$
\lim_{n \to \infty} P_\theta ( \theta' \in \hat{\Theta}_{1-\varepsilon}( Y^{(n)})) 
= 0\,.
$$

\vspace{.8cm}

\noindent
{\it Proof of (\ref{res})}.\\
 Using Laplace expansion (\ref{Laplace}) and the definition of $k_n$  we see that  
$$
\hat{\Theta}_{1-\varepsilon}( y^{(n)}) = \left\{ \theta \in \Theta \, : \; 2\left( \ell(\hat\theta_n;y^{(n)} ) - \ell(\theta;y^{(n)} ) \right) < 2 k_n +O(n^{-1})
 \right\}\,,  
$$
so that, if 
$$
2\left( \ell(\hat\theta_n;Y^{(n)} ) - \ell(\theta;Y^{(n)} ) \right) \stackrel{d}{\to} \chi^2_p\,,
$$
then (\ref{res}) follows from (\ref{kn}).

\subsection*{References}
\begin{description}

\item
$\,$
Barndorff--Nielsen, O.E.\ and Cox, D.R.\ (1989). {\it Asymptotic Techniques  for Use in Statistics}. London,  Chapman and Hall.

\item
$\,$
Benjamin, D.J., Berger, J.O., Johannesson, M., Nosek, B.A., Wagenmakers, E.J., Berk, R., et al.\  (2017). Redefine statistical significance. {\it Nature Human Behaviour},  {\bf 33}, 175. 

\item
$\,$
Bickel, D.R.\ (2019). Sharpen statistical significance: Evidence thresholds and Bayes factors sharpened into Occam's razor. {\it Stat}, {\bf 8}, e215.

\item
$\,$
Csenki, A.\ (1979). A note on confidence sequences in multiparameter exponential families.  
{\it Journal of Multivariate Analysis},  {\bf 9}, 337--340.

\item
$\,$
Darling, D.A.\ and Robbins, H. (1967a). Iterated logarithm inequalities. {\it Proceedings of the  National Academy of  Sciences of the  USA}, {\bf 57}, 1188--1192. 

\item
$\,$
Darling, D.A.\ and Robbins, H.\ (1967b). Confidence sequences for mean, variance and
median.  {\it Proceedings of the  National Academy of  Sciences of the  USA}, {\bf 58}, 66--68.

\item
$\,$
DiCiccio, T.J., Kuffner, T.A.\ and Young, G.A.\ (2015). Quantifying nuisance parameter effects via decompositions of asymptotic refinements for likelihood-based statistics. {\it Journal of Statistical Planning and Inference}, {\bf 165}, 1--12. 
%
%
%
%
%

\item
$\,$
Fay, M.P., Hunsberger, S.A., Nason, M. and Gabriel, E.\ (2018). exact2x2 -- exact tests and confidence intervals for 2x2 tables. {\it R Package Version 2018.07-27.} (Available from {\tt https://cran.r-project.org/web/packages/exact2x2/}.)

\item
$\,$
Fog, A.\ (2015). BiasedUrn -- biased urn model distributions. {\it R Package Version 2015.12-28.} (Available from {\tt https://cran.r-project.org/web/packages/BiasedUrn/}.)

\item
$\,$
Fraser, D.A.S., B\'{e}dard, M., Wong, A., Lin, W.\ and Fraser, A.M.\ (2016). Bayes, reproducibility and the quest for truth. {\it Statistical Science}, {\bf 31}, 578--590.

\item
$\,$
Gandy, A.\ and Hahn, G.\  (2016). A framework for Monte Carlo based multiple testing. {\it Scandinavian  Journal of  Statistics}, {\bf 43}, 1046--1063.  


\item
$\,$
Goodman, S.N.\ (2016). Aligning statistical and scientific reasoning. {\it Science}, {\bf 352}, 1180--1181. 

%
%
%
\item
$\,$
Ioannidis, J. P.A.\ (2005). Why most published research findings are false.  {\it PLoS Medicine},  {\bf 2}, e124. 


\item
$\,$
Jacod, J.\ and Protter, P.\ (2000). {\it Probability Essentials}. Berlin, Springer.  
\item
$\,$
Jennison, C.\ and Turnbull, B.W.\ (1989). Interim analyses: The repeated confidence interval approach.  {\it Journal of the Royal Statistical Society, Ser.\ B}, {\bf 51}, 305--361.  

\item
$\,$
Johari, R.J., Koomen, P., Pekelis, L.\ and Walsh, D.\ (2017). Peeking at A/B Tests: Why it matters, and what to do about it.  
In {\it Proceedings of the 23rd ACM SIGKDD International Conference on Knowledge Discovery and Data Mining}, 1517--1525, New York, ACM.

\item
$\,$
Johnson, V.E.\ (2013). Revised standards for statistical
evidence. {\it Proceedings of the  National Academy of  Sciences of the  USA}, {\bf 110}, 19313--19317.


\item
$\,$
Kuffner, T.A.\ and Walker, S.G.\  (2019). Why are p-values controversial? {\it The American Statistician}, {\bf 73}, 1--3.

\item
$\,$
Lai, T.L.\ (1976). On confidence sequences. {\it The Annals of  Statistics},  {\bf 4},  265--280.

\item $\,$
Lai, T.L.\ and Siegmund, D.\ (1986). The contributions of Herbert Robbins to mathematical statistics. {\it Statistical Science}, {\bf 1}, 276--284.

\item $\,$
 Lakens, D., Adolfi, F. G., Albers, C. J., Anvari, F., Apps, M.A.J., Argamon, S.E.,  van Assen, M.A.L. M.\ et al.\ (2017). Justify your alpha: A response to ``Redefine statistical significance''.  Retrieved from psyarxiv.com/9s3y6.

%
\item
$\,$
McCullagh, P.\ and Nelder, J.A.\ (1989). {\em Generalized Linear Models}, 
2-nd ed.. London, Chapman and Hall.

\item
$\,$
McShane, B.B.\ and Gal, D.\ (2017). Statistical significance and the dicho\-tom\-ization of evidence (with discussion). {\it Journal of the American Statistical Association}, {\bf 112}, 885--908. 

\item
$\,$
Morris, C.N.\ (1982). Natural exponential families with quadratic variance functions. {\it The Annals of Statistics}, {\bf 10}, 65--80. 
 
\item
$\,$
Nakagawa, S.\ and Cuthill, I.C.\ (2007). Effect size, confidence interval and statistical
significance: a practical guide for biologists. {\it Biological Reviews}, {\bf 82}, 591--605.


%

%
\item
$\,$
Robbins, H.\ (1970). Statistical methods related to the law of the iterated logarithm. {\it The Annals of Mathematical  Statistics},  {\bf 41},  1397--1409.

\item
$\,$
Royall, R.\  (1997). {\it Statistical Evidence: A Likelihood Paradigm}. London, 
Chapman and Hall.

%
\item
$\,$
Russel, S.J.\ and Norvig, P.\  (2010). {\it Artificial Intelligence: A Modern Approach}, Third Ed., Prentice Hall, NJ.

\item
$\,$
Schwalbe, M.\ (2016). {\it Statistical Challenges in Assessing and Fostering the Reproducibility of Scientific Results:
Summary of a Workshop}. Washington, DC: National Academies Press.

\item
$\,$
Severini, T.A.\ (2000). {\it Likelihood Methods in Statistics}, Oxford University Press, Oxford.

\item
$\,$
Wasserstein, R.L.\ and  Lazar, N.A.\  (2016). The ASA's statement on
p-values: Context, process, and purpose. {\it The American Statistician}, {\bf 70}, 129--133.

\item
$\,$
Whitehead, J.\ (1989). Discussion of:   Interim analyses: The repeated confidence interval approach by C.\ Jennison and B.W.\ Turnbull,  {\it Journal of the Royal Statistical Society, Ser.\ B}, {\bf 51}, 338.  

\end{description}

%
%


\end{document}